\documentclass[12pt]{amsart}

\usepackage{amsthm, amsmath, amssymb}
\usepackage{comment}
\usepackage{color}

\newtheorem{theorem}{Theorem}
\newtheorem{lemma}{Lemma}
\newtheorem{proposition}{Proposition}
\newtheorem{corollary}{Corollary}
\theoremstyle{definition}
\newtheorem{definition}{Definition}

\theoremstyle{remark}
\newtheorem{remark}{Remark}

\date{}

\author{Avetik Arakelyan}
\address{Institute of Mathematics, NAS of Armenia, 0019 Yerevan, Armenia}
\email{arakelyanavetik@gmail.com}

\thanks{A. Arakelyan was supported by State Committee of Science MES RA, in frame of the research project No.  16YR-1A017 }

\title[ Avetik Arakelyan]{Convergence of  the Finite  Difference Scheme for a General Class of the Spatial Segregation of Reaction-diffusion Systems}



\begin{document}

\begin{abstract}
	In this work we prove convergence of the finite difference scheme for equations of stationary states of a general class of the spatial segregation of reaction-diffusion systems with $m\geq 2$ components. More precisely, we show  that  the numerical solution $u_h^l$, given by the difference scheme, converges to the $l^{th}$ component $u_l,$ when  the mesh size $h$ tends to zero, provided $u_l\in C^2(\Omega),$ for every $l=1,2,\dots,m.$  In particular, our proof provides convergence of a difference scheme for the  multi-phase obstacle problem. 
\end{abstract}

\maketitle

%\tableofcontents

\section{Introduction}

\subsection{The  setting of the problem }
In recent years there have been intense studies of spatial segregation for reaction-diffusion systems. The existence of spatially inhomogeneous solutions for competition models of Lotka-Volterra type in the case of two and more competing densities  have been considered in
\cite{MR2146353,MR2151234,MR2283921,MR2300320,MR1900331,MR1459414, MR2417905}.
Aforementioned segregation problems led to an interesting class of  multi-phase obstacle-like free boundary problems. These   problems have growing interest due to their  important applications in the different branches of applied mathematics. To see the diversity of applications we refer \cite{Avet-Henrik,bucur-multi,Aram-Caff} and the references therein.  

Nowadays, the theory of the one- and  two-phase obstacle-like problems (elliptic and parabolic versions) is well-established  and for a reference we address to the books \cite{PSU2012,rodrigues-book} and references therein. For two-phase problems the interested reader is also referred to the recent works \cite{faridperturb,repin2015}.

 There is a vast literature devoted to the  numerical analysis of one-phase obstacle-like problems, and we refer some of well-known papers \cite{nochetto-residual,multigrid,American-valuationl,oberman-siam}.   For the numerical treatment of the two-phase  problems we refer to the works \cite{arakelyan2015finite,ABP2014,MR2961456,Farid,farid-fem,Osher}.

 The present work concerns to  prove the convergence of the difference scheme for a certain class of the spatial segregation of reaction-diffusion system with $m$ components.

Let $\Omega  \subset \mathbb{R}^n, n\geq 2$ be a connected and  bounded domain with smooth boundary and  $m$ be a fixed integer.
We consider the steady-states of $m$ competing species coexisting in  the same area $\Omega$.
Let $u_{i}(x)$ denotes  the population density of the  $i^\textrm{th}$ component with the internal dynamic  prescribed by $F_{i}(x,u_i)$.

We call the $m$-tuple $U=(u_1,\cdots,u_m)\in (W^{1,2}(\Omega))^{m},$ \emph{a segregated state} if
\[
u_{i}(x) \cdot  u_{j}(x)=0,\  \text{a.e. } \text{ for  } \quad i\neq j,  \ x\in \Omega.
\]
The problem amounts to
\begin{equation}\label{main_problem}
\text{  Minimize  }  E(u_1, \cdots, u_m)=\int_{\Omega}  \sum_{i=1}^{m} \left( \frac{1}{2}| \nabla u_{i}|^{2}+F_i(x,u_i) \right) dx,
\end{equation}
  over the set
  $$S={\{(u_1,\dots,u_{m})\in (W^{1,2}   (\Omega))^{m} :u_{i}\geq0, u_{i} \cdot u_{j}=0, u_{i}=\phi_{i} \quad \text {on} \quad \partial  \Omega}\},$$
where $\phi_{i} \in H^{\frac{1}{2}}(\partial \Omega),$\; $\phi_{i}  \cdot \phi_{j}=0,$ for $i\neq j$ and $ \phi_{i}\geq 0$ on the boundary $\partial \Omega.$

We assume that
\[F_i(x,s)=\int_0^s f_i(x,v)dv,\]
 where $f_i(x,s):\Omega\times\mathbb{R}^+\to\mathbb{R}$ is Lipschitz continuous in $s,$ uniformly continuous in $x$ and $f_i(x,0)\equiv0.$
%In the sequel, we assume that the functional  \eqref{main_problem} is coercive, which will be needed to provide the existence of minimizers. In order to have  coercivity for the functional, for instance following \cite{ MR2151234}, one can assume that
%there exists $b_i(x)\in L_\infty(\Omega)$ such that
%\[
%|f_i(x,s)|\leq b_i(x)\cdot s, \;\;\forall x\in\Omega,\;\; s\geq s_0 >> 1,
%\]
%and
%\[
%\int_\Omega(|\nabla w(x)|^2 - b_i(x)w^2(x))dx>0, \quad\forall w \in H_0^1(\Omega)\setminus\{0\}.
%\]

\begin{remark}
Functions $f_i(x,s)$'s are defined only for non negative values of s (recall that our densities $u_i$'s are assumed non negative); thus we can arbitrarily define such functions on the negative semiaxis. For the sake of convenience, when $s\leq 0$,  we will let $f_i(x,s)=-f_i(x,-s).$ This extension  preserves the continuity due to the conditions on $f_i$ defined above. In the same way, each $F_i$ is extended as an even function.
\end{remark}
\begin{remark}
We emphasize that for the case $f_i(x,s)=f_i(x),$ the  assumption is that  for all $i$ the functions $f_i(x,s)$ are nonnegative and  uniformly continuous in $x.$  Also for simplicity, throughout the paper we shall call both $F_i(x,u_i)$ and $f_i(x,u_i)$ internal dynamics.

%And throughout the paper we will  consider the case $f_i(x,s)$ depends on the variable $s,$ keeping in mind that the same results will take place also when  $f_i(x,s)$  doesn't depend on $s.$ 
\end{remark}

	We would like to point out that the only difference between our minimization problem \eqref{main_problem} and the problem discussed in \cite{MR2151234}, is  the sign in front of the internal dynamics $F_i$.  In our case, the plus sign of $F_i$ allows to get rid of some additional  conditions, which are imposed in \cite[Section $2$]{MR2151234}. Those conditions are important to provide coercivity of a minimizing functional in \cite{MR2151234}. But in our case the above given conditions together with convexity assumption on $F_i(x,s),$ with respect to the variable $s$ are enough to conclude $F_i(x,u_i(x))\geq 0,$ which in turn implies coercivity of a functional \eqref{main_problem}.

In order to speak on  the local properties of the population densities, let us introduce the notion of multiplicity of a point in $\Omega$.
\begin{definition}
	The multiplicity of the point $x \in \overline{\Omega}$ is defined by:
	\begin{equation*}
		m(x)=\text{card} \left\{i: measure (\Omega_{i} \cap B(x,r))>0, \forall r>0\right\},
	\end{equation*}
	where  $\Omega_i=\{u_i>0\}$.
\end{definition}

For the local properties of $u_i$ the same results as in \cite{MR2151234} with  the opposite sign in front of the internal dynamics $f_i$ hold. Below, for the sake of clarity, we  write down these results from \cite{MR2151234} with appropriate changes.

\begin{lemma}(Proposition $6.3$ in \cite{MR2151234})\label{CTV-thm1}
Assume that $x_{0} \in \Omega,$ then the following holds:\\
\begin{enumerate}
\item[1)] If $m(x_0)=0,$ then there exists $r>0$ such that for every $i=1,\cdots m$;
\[
 u_i\equiv0\;\; \mbox{on}\;\;  B(x_0,r).
\]
\item[2)] If $m(x_0)=1,$ then there are $i$  and $r>0$ such that in $ B(x_0,r)$
\[
\Delta u_i=f_i(x,u_i), \quad \ \quad u_j\equiv0 \quad \text{for  } j\neq i.
\]
\item[3)] If $m(x_0)=2,$ then there are $i, j $ and  $r>0$ such that for every $k$  and $k\neq i,j$, we have $ u_k\equiv 0$ and
\[
\Delta(u_i-u_j) =f_{i}(x,(u_i-u_j))\chi_{\{u_i>u_j\}} -   f_{j}(x,-(u_i-u_j))\chi_{\{u_i<u_j\}}  \text{ in  }  B(x_0,r).
\]
\end{enumerate}
\end{lemma}
\begin{lemma}(Theorem $5.1$ in \cite{ MR2151234})\label{CTV-thm2}
For every minimizer $(u_1,\cdots,u_m)\in S$ to the functional  \eqref{main_problem}, the following inequality holds
$$
\Delta\left(u_l(x)-\sum_{p\neq l}u_p(x)\right)\leq f_l(x,u_l),
$$
for all $l=1,2,\dots,m.$

\end{lemma}

Next, we state the following uniqueness theorem due to  Conti, Terrachini and Verzini, by observing that in our case  the plus sign in front of $F_i$  requires convexity condition on $F_i(x,s)$ rather than concavity condition given in \cite{MR2151234}.
\begin{theorem}(Theorem $4.2$ in \cite{ MR2151234})
Let the functional in minimization problem \eqref{main_problem} is coercive and moreover each $F_i(x,s)$ is convex in the variable $s,$ for all $x\in\Omega.$ Then, the problem  \eqref{main_problem} has a unique minimizer.
\end{theorem}

\subsection{Notation}

We will work in two-dimensional space $\mathbb{R}^2.$ For the sake of simplicity, we will assume that $\Omega=(0,a)\times(0,a).$  It should be remarked that the same results can be obtained rigorously also for more complicated domains. 

Let $N\in\mathbb{N}$ be a positive integer, $h=a/N$ and
$$
x_i=ih,\, y_i=ih,\quad i=0,1,...,N.
$$

We use the notation $u^l_h(x,y)$ for  the finite difference scheme approximation to  $u_l(x,y).$ We will heavily use the shorthand notations $u^l_h(z)$ and $u_l(z),$ where $z=(x,y)\in\Omega.$ Concerning the boundary functions $\phi_l,$ we assume they are  extended to be zero everywhere outside the boundary $\partial\Omega,$ for all $l=1,2,\dots, m.$ The discrete approximation for these functions will be  $\phi^l_h.$

%In this paper we will also use  notations $u=(u_{\alpha})$,  (not to be confused with functions $u$).% Also we will write $(a_{\alpha})\le (b_{\alpha})$ in $\mathcal I$ if $a_{\alpha}\le b_{\alpha}$ for all $\alpha\in\mathcal I$.

Denote
$$
 \Omega_h=\{(x_i,y_j)=(ih,jh):\ 0\leq i,j\leq N\},
$$
$$
   \Omega_h^o=\{(x_i,y_j)=(ih,jh):\ 1\leq i,j\leq N-1\},
$$
and
$$
  \partial \Omega_h= \Omega_h \setminus  \Omega_h^o.
$$

In  two-dimensional case we introduce the following 5-point stencil approximation for Laplacian:
$$
 L_h v(x,y)=\frac{v(x-h,y)+v(x+h,y)-4v(x,y)+v(x,y-h)+v(x,y+h)}{h^2}
$$
for any $(x,y)\in \Omega$.

\section{Finite difference scheme }
We start this section by defining the finite difference scheme, which convergence analysis will be the subject of the study in the present work. We denote it by $(u_h^1,u_h^2,\dots,u_h^m).$ This vector solves the following system:
\begin{equation}\label{scheme_sys}
\begin{cases}
u_h^{l}(z) =\max \left(\cfrac{-f_{l}(z,u_h^{l}(z))h^{2}}{4}+\overline{u}_h^l(z) - \sum\limits_{p \neq l}  \overline{u}_h^p(z) , \,  0\right),\;\;z\in\Omega_h^o,\\
u_h^{l}(z) =\phi_h^{l}(z)= \phi_{l}(z),\;\;z\in \partial\Omega_h.
\end{cases}
\end{equation}
for every $l=1,2,\dots,m$ and $z=(x,y)\in\Omega_h.$ 
Here for a given uniform  mesh on $\Omega\subset \mathbb{R}^2,$   we define     $\overline{u}_h^l(z)$ to be the average of $u_h^{l}(z)$ for all neighbor points of $z=(x,y)\in\Omega_h^o:$
\[
\overline{u}_h^l(z)=\frac{1}{4}[u_h^{l}(x_{i-1},y_{j})+u_h^{l} (x_{i+1},y_{j})+u_h^{l}(x_{i},y_{j-1})  +u_h^{l} (x_{i},y_{j+1})].
\]
Throughout the paper the following notations will play a crucial role:
$$
\hat{u}_{l}(z):={u}_{l}(z)-\sum_{p\neq l} {u}_{p}(z),
$$
and
$$
\hat{u}_h^{l}(z):={u}_h^{l}(z)-\sum_{p\neq l} {u}_h^{p}(z).
$$
It is easy to verify that the solution $({u}_h^{1},{u}_h^{1},\dots,{u}_h^{m})$ to a difference scheme \eqref{scheme_sys} for every $l=1,2,\dots,m,$ satisfies the following properties, provided that all functions $f_{l}(z,s)$ are nondecreasing with respect to the variable $s:$
\begin{align}\label{scheme-prop}
	\begin{cases}
	   L_h(\hat{u}_h^{l}(z))\leq f_{l}(z,u_h^{l}(z))\;\;&z\in\Omega_h,\\
	   L_h(\hat{u}_h^{l}(z))= f_{l}(z,u_h^{l}(z))\;\;&z\in\{u_h^{l}(z)>0\},\\
	   {u}_h^{l}(z)\geq 0\;\;&z\in\Omega_h,\\
	   {u}_h^{p}(z)\cdot {u}_h^{q}(z)=0,\;\;p\neq q,\\
	   {u}_h^{l}(z)={\phi}_h^{l}(z)={\phi}_{l}(z),\;\;&z\in\partial\Omega_h.
	\end{cases}
\end{align}

  The  difference system \eqref{scheme_sys}, when the internal dynamics  $f_i(z,s)=0,\;z\equiv(x,y)\in \mathbb{R}^2,$ has been suggested in \cite{MR2563520}.  The author only implemented plausible numerical  figures by this scheme, without  its analysis. This finite difference method has been generalized in \cite{Mywork} for the case of non-negative internal dynamics $f_i(z,s)=f_i(z).$ In \cite{Mywork} the authors give a numerical consistent variational system with strong interaction, and provide disjointness condition of populations during the iteration of the scheme. In this case the proposed algorithm is lack of deep analysis,  especially for the case of three and more competing populations. In the recent work by the current author in collaboration \cite{Av-Raf-2016} the existence and uniqueness of the scheme, which solves the system \eqref{scheme_sys}, have been proven, provided all $f_i(z,s)$ are nonnegative and nondecreasing with respect to $s.$ It is noteworthy, that the difference schemes with the same spirit as the system \eqref{scheme_sys}, have been successfully applied in  quadrature domains theory (see \cite{farid-multiQD}) and in optimal partitions theory (see \cite{farid-optimal-partition}). This makes us to strongly believe that the ideas behind the difference scheme \eqref{scheme_sys} have great opportunities  to be applied in different problems, where the segregated geometry arise.

\section{Auxiliary lemmas}\label{unique_schm}
In this section we prove two technical lemmas, which will be used for the convergence analysis of the scheme. To this aim, for the sake of convenience  we  denote by $nbr(z)$ the set of all closest neighbor points corresponding to a mesh point $z=(x,y)\in\Omega_h.$ We will need also the following barrier function:
$$
V_h(z)=V_h(x,y)=\cfrac{1}{4}(x^2+y^2+1)\sum\limits_{l=1}^{m}||L_hu_l-\Delta u_l||_{L^\infty(\Omega)}.
$$  
For simplicity, we set  by
$$
\alpha_h:=L_h(V_h(z))=\sum\limits_{l=1}^{m}||L_hu_l-\Delta u_l||_{L^\infty(\Omega)}.
$$
\begin{lemma}\label{lemma1}
	Let the functions $f_{l}(z,s)$ be nondecreasing with respect to the variable $s.$ We set $(u_{1},u_{2},\dots,u_{m})\in S\cap (C^2(\Omega))^m$ to be an exact minimizer of \eqref{main_problem} subject to $S,$  and by  $(u_h^{1},u_h^{2},\dots,u_h^{m})$ we define the  vector, which solves the finite difference system \eqref{scheme_sys}. Then the following statements are true:
	\[
	\max_{\Omega_h}\left(\hat{u}_l(z)-\hat{u}_h^l(z)+V_h(z)\right)=\max_{\{{u}_l(z)\leq u_h^l(z)\}}\left(\hat{u}_l(z)-\hat{u}_h^l(z)+V_h(z)\right),
	\]
	and
	\[
	\max_{\Omega_h}\left( \hat{u}_h^l(z)-\hat{u}_l(z)+V_h(z)\right)=\max_{\{{u}_h^l(z)\leq u_l(z)\}}\left( \hat{u}_h^l(z)-\hat{u}_l(z)+V_h(z)\right),
	\]
	for all $l=1,2,\dots,m$.
\end{lemma}

\begin{proof}
	We argue by contradiction. Suppose for some $l_0$ we have
	\begin{multline}\label{init_assmp}
	\hat{u}_{l_0}(z_0)-\hat{u}_h^{l_0}(z_0)+V_h(z_0)=
	\max_{\Omega_h}\left(\hat{u}_{l_0}(z)-\hat{u}_h^{l_0}(z)+V_h(z)\right)=\\=\max_{\{{u}_{l_0}(z)> u_h^{l_0}(z)\}}\left(\hat{u}_{l_0}(z)-\hat{u}_h^{l_0}(z)+V_h(z)\right)>\\>
	\max_{\{{u}_{l_0}(z)\leq u_h^{l_0}(z)\}}\left(\hat{u}_{l_0}(z)-\hat{u}_h^{l_0}(z)+V_h(z)\right).
	\end{multline}
	Then taking into account the following simple chain of inclusions
	\begin{equation}\label{incl_chain}
		\{{u}_l(z)> u_h^l(z)\}\subset\{\hat{u}_l(z)> \hat{u}_h^l(z)\}\subset\{{u}_l(z)\geq u_h^l(z)\},
	\end{equation}
	we obviously see that $ {u}_{l_0}(z_0)> u_h^{l_0}(z_0)\geq 0 $ implies 
	$\hat{u}_{l_0}(z_0)> \hat{u}_h^{l_0}(z_0).$ On the other hand,  the discrete system \eqref{scheme_sys} and Lemma \ref{CTV-thm1} gives us
	 $$
	\Delta\hat{u}_{l_0}(z_0)=\Delta{u}_{l_0}(z_0)=f_{l_0}(z_{0},{u}_{l_0}(z_0))\;\; \mbox{and}\;\; L_h\hat{u}_h^{l_0}(z_0)\leq f_{l_0}(z_{0},{u}_h^{l_0}(z_0)).$$ 
	Therefore 
	\begin{multline*}
		L_h\left(\hat{u}_{l_0}(z_0)-\hat{u}_h^{l_0}(z_0)+V_h(z_0)\right)=\\= \Big(L_h\hat{u}_{l_0}(z_0)-\Delta\hat{u}_{l_0}(z_0)+\alpha_h \Big)+\left(\Delta\hat{u}_{l_0}(z_0)- L_h\hat{u}_h^{l_0}(z_0)\right)\geq\\\geq \left(\Delta\hat{u}_{l_0}(z_0)- L_h\hat{u}_h^{l_0}(z_0)\right)\geq f_{l_0}(z_0,u_{l_0}(z_0))-f_{l_0}(z_0,u_h^{l_0}(z_0)) \geq 0.
	\end{multline*}
	Thus,
	$$
	\hat{u}_{l_0}(z_0)-\hat{u}_h^{l_0}(z_0)+V_h(z_0)\leq \frac{1}{4}\sum_{\{z\in nbr(z_0)\}}\left(\hat{u}_{l_0}(z)-\hat{u}_h^{l_0}(z)+V_h(z)\right),
	$$
	which implies that $\hat{u}_{l_0}(z_0)-\hat{u}_h^{l_0}(z_0)+V_h(z_0)=\hat{u}_{l_0}(z)-\hat{u}_h^{l_0}(z)+V_h(z),$ for all $z\in nbr(z_0).$ Since $\hat{u}_{l_0}(z_0)> \hat{u}_h^{l_0}(z_0),$ then we apparently have  
	$$
	\hat{u}_{l_0}(z)-\hat{u}_h^{l_0}(z)>V_h(z_0)-V_h(z),
	$$ 
	for all ${z}\in nbr(z_0).$  We take a particular neighbor point $\acute{z}=(x_{{i_0}-1},y_{j_0}),$ provided $z_0=(x_{i_0},y_{j_0})\in\Omega_h.$ 
	We obtain
	$$
	\hat{u}_{l_0}(\acute{z})-\hat{u}_h^{l_0}(\acute{z})>V_h(z_0)-V_h(\acute{z})=
	\frac{1}{4}(x_{i_0}^2-x_{{i_0}-1}^2)\alpha_h\geq 0.
	$$
	In view of  chain \eqref{incl_chain} we get ${u}_{l_0}(\acute{z})\geq u_h^{l_0}(\acute{z}).$ 
	According to our assumption \eqref{init_assmp}, the only possibility is ${u}_{l_0}(\acute{z})> u_h^{l_0}(\acute{z}).$  Now we can proceed the previous steps for this neighbor point $\acute{z}=(x_{{i_0}-1},y_{j_0})\in nbr(z_0),$ and obtain  the same strict inequality for $(x_{{i_0}-2},y_{j_0})$ and so on. Continuing this along an $x$ axis, we will finally approach to the boundary $\partial\Omega_h,$ where as we know ${u}_{l_0}(z)={u}_h^{l_0}(z)={\phi}_{l_0}(z),$ for all $z\in\partial\Omega_h.$ Hence, the strict inequality  fails, which implies that our initial assumption \eqref{init_assmp} is false. Observe that the same arguments can be applied if we interchange the role of ${u}_{l}(z)$ and ${u}_h^{l}(z).$ In this case we need to use  the reversed chain of inclusions given below
	\begin{equation}\label{chain_reverse}
	\{{u}_h^l(z)> u_l(z)\}\subset\{\hat{u}_h^l(z)> \hat{u}_l(z)\}\subset\{{u}_h^l(z)\geq u_l(z)\},
	\end{equation}
	and  Lemma \ref{CTV-thm2}. Thus, we also have 
	\[
	\max_{\Omega_h}\left( \hat{u}_h^l(z)-\hat{u}_l(z)+V_h(z)\right)=\max_{\{{u}_h^l(z)\leq u_l(z)\}}\left( \hat{u}_h^l(z)-\hat{u}_l(z)+V_h(z)\right),
	\]
	for every $l=1,2,\dots,m.$ 	This completes the proof of Lemma.
\end{proof}
%	It is easy to see that $\max\limits_{\mathcal N}(\hat{v}_\alpha^l-\hat{u}_\alpha^l)=-\min\limits_{\mathcal N}(\hat{u}_\alpha^l-\hat{v}_\alpha^l).$ 
%	Particularly, for every fixed $l=1,2.\dots,m$ and $\alpha\in\mathcal N$ we have
%	\begin{equation}\label{double_ineq}
%	-\max_{\{v_\alpha^l\leq u_\alpha^l\}}(\hat{v}_\alpha^l-\hat{u}_\alpha^l)=
%	-\max\limits_{\mathcal N}(\hat{v}_\alpha^l-\hat{u}_\alpha^l)\leq \hat{u}_\alpha^l-\hat{v}_\alpha^l\leq  \max\limits_{\mathcal N}(\hat{u}_\alpha^l-\hat{v}_\alpha^l)=\max_{\{u_\alpha^l\leq v_\alpha^l\}}(\hat{u}_\alpha^l-\hat{v}_\alpha^l).
%	\end{equation}
In the sequel and thanks to Lemma \ref{lemma1}, we will use the following notations:

\begin{align*}
M_h:&=	
\max_l\left(\max_{\Omega_h}\left(\hat{u}_l(z)-\hat{u}_h^l(z)+V_h(z)\right)\right)\\&=\max_l\left(\max_{\{{u}_l(z)\leq u_h^l(z)\}}\left(\hat{u}_l(z)-\hat{u}_h^l(z)+V_h(z)\right)\right),
\end{align*}
and
\begin{align*}
R_h:&=\max_l\left(\max_{\Omega_h}\left( \hat{u}_h^l(z)-\hat{u}_l(z)+V_h(z)\right)
\right)\\&=\max_l\left(\max_{\{{u}_h^l(z)\leq u_l(z)\}}\left( \hat{u}_h^l(z)-\hat{u}_l(z)+V_h(z)\right)
\right).
\end{align*}

\begin{lemma}\label{lemma2}
	Let the functions $f_{l}(x,s)$ be nondecreasing with respect to the variable $s.$ We also set  $(u_{1},u_{2},\dots,u_{m})\in S\cap (C^2(\Omega))^m$ to be an exact minimizer of \eqref{main_problem} subject to $S,$ and $(u_h^{1},u_h^{2},\dots,u_h^{m})$  to be the difference scheme, which solves  the discrete system \eqref{scheme_sys}. For these two elements we set $M_h$ and $R_h$ as defined above. If $M_h>\max_{\Omega_h}V_h(z)$(respectively $R_h>\max_{\Omega_h}V_h(z)$) and it is attained for some $l_0$, then $M_h=R_h>\max_{\Omega_h}V_h(z).$ Moreover, there exists some $t_0\neq l_0,$ and $z_0\in\Omega_h,$ such that
$$
M_h=\max_{\{{u}_{l_0}(z)= u_h^{l_0}(z)=0\}}\left(\hat{u}_{l_0}(z)-\hat{u}_h^{l_0}(z)+V_h(z)\right)={u}_h^{t_0}(z_0)-{u}_{t_0}(z_0)+V_h(z_0).
$$
(Respectively, 
$$
	R_h=\max_{\{{u}_{l_0}(z)= u_h^{l_0}(z)=0\}}\left(\hat{u}_h^{l_0}(z)-\hat{u}_{l_0}(z)+V_h(z)
	\right)={u}_{t_0}(z_0)-{u}_h^{t_0}(z_0)+V_h(z_0)).
$$
\end{lemma}

\begin{proof}
	Due to Lemma \ref{lemma1} we have 
$$
M_h=\max_{\{{u}_{l_0}(z)\leq u_h^{l_0}(z)\}}\left(\hat{u}_{l_0}(z)-\hat{u}_h^{l_0}(z)+V_h(z)\right).
$$

It is easy to verify that $(\hat{u}_{l_0}(z)-\hat{u}_h^{l_0}(z))$ might be strictly positive only on the set $\{{u}_{l_0}(z)= u_h^{l_0}(z)=0\}$ (for the other cases $(\hat{u}_{l_0}(z)-\hat{u}_h^{l_0}(z))\leq 0).$ Hence, 
	$$
      \left(\hat{u}_{l_0}(z)-\hat{u}_h^{l_0}(z)+V_h(z)\right)\leq \max_{\Omega_h}V_h(z)< M_h,
	$$
	provided $z\notin \{{u}_{l_0}(z)= u_h^{l_0}(z)=0\},$ which yields 
	$$
	 M_h=\max_{\{{u}_{l_0}(z)= u_h^{l_0}(z)=0\}}\left(\hat{u}_{l_0}(z)-\hat{u}_h^{l_0}(z)+V_h(z)\right).
	 $$
	Using the latter equality, one can prove that $M_h>\max_{\Omega_h}V_h(z)$ implies $R_h>\max_{\Omega_h}V_h(z).$  Indeed, it is easy to see that if the maximum $M_h$ is attained at the mesh point $z_0\in\Omega_h,$ then there exists  $t_0\neq l_0$ such that
	\begin{multline}\label{M=R}
	\max_{\Omega_h}V_h(z)< M_h=\max_{\{{u}_{l_0}(z)= u_h^{l_0}(z)=0\}}\left(\hat{u}_{l_0}(z)-\hat{u}_h^{l_0}(z)+V_h(z)\right)=\\=
	\hat{u}_{l_0}(z_0)-\hat{u}_h^{l_0}(z_0)+V_h(z_0)=
	\sum_{l\neq l_0}\left({u}_h^{l}(z_0)-{u}_{l}(z_0)\right)+V_h(z_0)
	=\\={u}_h^{t_0}(z_0)-\sum_{l\neq l_0}{u}_{l}(z_0)+V_h(z_0)\leq  \hat{u}_h^{t_0}(z_0)-\hat{u}_{t_0}(z_0)+V_h(z_0)\leq R_h.
	\end{multline}
	In the same way we will obtain that $\max_{\Omega_h}V_h(z)<R_h\leq M_h,$ and therefore $$M_h=R_h>\max_{\Omega_h}V_h(z).$$
	On the other hand, the above computation \eqref{M=R} gives us
	
	$$
	{u}_h^{t_0}(z_0)-\sum_{l\neq l_0}{u}_{l}(z_0)
	=\hat{u}_h^{t_0}(z_0)-\hat{u}_{t_0}(z_0).
	$$ 
	This leads to $2\sum\limits_{l\neq t_0}{u}_{l}(z_0)= 0,$ and therefore
	${u}_{l}(z_0)=0,$ for all $l\neq t_0.$ Hence,
	$$M_h={u}_h^{t_0}(z_0)-\sum_{l\neq l_0}{u}_{l}(z_0)+V_h(z_0)={u}_h^{t_0}(z_0)-{u}_{t_0}(z_0)+V_h(z_0).
	$$
	For $R_h$ the proof can be done in a similar way.
	This completes the proof.
\end{proof}

\begin{comment}

Lets consider the particular case  $m=2,$ when the system \eqref{scheme_sys}  reduces to
 \begin{equation}\label{scheme_sys_two}
 \begin{cases}
 w_\alpha^{1} =\max \left(\cfrac{-f_{1}(x_\alpha,w_\alpha^{1})h^{2}}{4}+\overline{w}_\alpha^{1} - \overline{w}_\alpha^{2} , \,  0\right),\;\;\alpha\in\mathcal N^o,\\
  w_\alpha^{2} =\max \left(\cfrac{-f_{2}(x_\alpha,w_\alpha^{2})h^{2}}{4}+\overline{w}_\alpha^{2} - \overline{w}_\alpha^{1} , \,  0\right),\;\;\alpha\in\mathcal N^o,\\
 w_\alpha^{1}=\phi_\alpha^1,\;\;\alpha\in\mathcal \partial\mathcal N,\\
 w_\alpha^{2}=\phi_\alpha^2,\;\;\alpha\in\mathcal \partial\mathcal N,
 \end{cases}
 \end{equation}
 for every  $\alpha\in\mathcal N.$ 
 
It is noteworthy that  this system, where $f_{l}(x,s)=f_{l}(x),\; l=1,2,$  has already been suggested in \cite{Farid} to approximate the solution of a two-phase obstacle (membrane) problem. The convergence of the iterative algorithm corresponding to this particular case of \eqref{scheme_sys_two} was proved in \cite{avetikconv}. 

In view of Section  \ref{twopop_section} we clearly see  two equivalent ways to define a difference method for the segregation problem with two competing densities. The  first  approach is to consider the solution to a discrete Min-Max equation \eqref{nonlinprob}, and the second approach is the solution to a discrete system \eqref{scheme_sys_two}. 
\end{comment}	
\section{Convergence of scheme}

In this section we prove the main result of the paper.  Next proposition shows the estimate between the exact and numerical solutions. Then the pointwise  convergence of the scheme follows immediately. 

\begin{proposition}\label{scheme_conv}
	Let the functions $f_{l}(x,s)$ be nondecreasing with respect to the variable $s.$ We set $(u_{1},u_{2},\dots,u_{m})\in S$ to be an exact minimizer of \eqref{main_problem} subject to $S.$ If  $(u_h^{1},u_h^{2},\dots,u_h^{m})$ is the difference scheme, which solves  the discrete system \eqref{scheme_sys},  then the following estimate holds: 
	$$
	||u_l-u_h^l||_{L^\infty(\Omega_h)}\leq C_{\Omega}\cdot\sum_{l=1}^{m}||L_h u_l-\Delta u_l||_{L^\infty(\Omega)},
	$$ 
	for every $l\in\{1,2,\dots,m\},$ provided $(u_{1},u_{2},\dots,u_{m})\in S\cap (C^2(\Omega))^m.$ Here $C_{\Omega}>0$ is a constant depending only on $\Omega.$ 
\end{proposition}

\begin{proof}
  	
  For the  vectors $(u_{1},u_{2},\dots,u_{m})$ and $(u_h^{1},u_h^{2},\dots,u_h^{m})$ we set the definition of $M_h$ and $R_h.$ We are going to prove that $M_h\leq \max_{\Omega_h}V_h(z).$ As a consequence we will obtain that $R_h\leq \max_{\Omega_h}V_h(z)$ holds as well.

	Suppose $M_h> \max_{\Omega_h}V_h(z).$ Our aim is to prove that this case leads to a contradiction. Let the value $M_h$ is attained for some $l_0\in\overline{1,m},$ then
	due to Lemma \ref{lemma2}  we  have $M_h=R_h,$ and there exist $z_0\in \Omega_h$ and $t_0\neq l_0$ such that:
	\begin{align*}
		M_h=R_h&=\max_{\{{u}_{l_0}(z)u_h^{l_0}(z)=0\}}\left(\hat{u}_{l_0}(z)-\hat{u}_h^{l_0}(z)+V_h(z)\right)\\&={u}_h^{t_0}(z_0)-{u}_{t_0}(z_0)+V_h(z_0).
	\end{align*}

	This yields
	$$
	{u}_h^{t_0}(z_0)-{u}_{t_0}(z_0)=M_h-V_h(z_0)>\max_{\Omega_h}V_h(z)-V_h(z_0)\geq 0,
	$$
	and therefore  due to \eqref{scheme-prop} and Lemma \ref{CTV-thm2} we clearly obtain 
	$$ L_h\hat{u}_h^{t_0}(z_0)=f_{t_0}(z_{0},{u}_h^{t_0}(z_0))\;\;\mbox{and}\;\; 
	\Delta\hat{u}_{t_0}(z_0)\leq f_{t_0}(z_{0},{u}_{t_0}(z_0)).$$  
	By proceeding similar steps as in the proof of Lemma \ref{lemma1} and recalling that $f_{l}(x,s)$ are nondecreasing with respect to the variable $s$,  we conclude
	$$
	L_h\left(\hat{u}_h^{t_0}(z_0)-\hat{u}_{t_0}(z_0)+V_h(z_0)\right)\geq 0.
	$$
	Thus, 
	$$
	\hat{u}_h^{t_0}(z_0)-\hat{u}_{t_0}(z_0)+V_h(z_0)\leq \frac{1}{4}\sum_{\{\gamma\in nbr(z_0)\}}(\hat{u}_h^{t_0}(\gamma)-\hat{u}_{t_0}(\gamma)+V_h(\gamma)),
	$$
	which implies that $$M_h=\hat{u}_h^{t_0}(z_0)-\hat{u}_{t_0}(z_0)+V_h(z_0)=\hat{u}_h^{t_0}(\gamma)-\hat{u}_{t_0}(\gamma)+V_h(\gamma)>\max_{\Omega_h}V_h(z),$$
 for all $\gamma\in nbr(z_0).$ Hence, $\hat{u}_h^{t_0}(\gamma)>\hat{u}_{t_0}(\gamma)$ and this along with the chain \eqref{chain_reverse} gives that  for all $\gamma\in nbr(z_0),$ we have ${u}_h^{t_0}(\gamma)\geq{u}_{t_0}(\gamma).$ For the neighbor mesh points $\gamma$ we proceed as follows: If  ${u}_h^{t_0}(\gamma)>{u}_{t_0}(\gamma),$ for some $\gamma_0\in nbr(z_0),$ then obviously $$L_h\left(\hat{u}_h^{t_0}(\gamma_0)-\hat{u}_{t_0}(\gamma_0)+V_h(\gamma_0)\right)\geq 0.$$ 
	This, as we saw a few lines above, leads to 
	\begin{equation}\label{nbr-cond1}
	M_h=\hat{u}_h^{t_0}(\gamma_0)-\hat{u}_{t_0}(\gamma_0)+V_h(\gamma_0)=\hat{u}_h^{t_0}(\theta)-\hat{u}_{t_0}(\theta)+V_h(\theta)>\max_{\Omega_h}V_h(z),
	\end{equation}
	for all $\theta\in nbr(\gamma_0).$

	If   ${u}_h^{t_0}(\gamma)={u}_{t_0}(\gamma),$ for some $\gamma_0\in nbr(z_0),$ then due to $$\hat{u}_h^{t_0}(\gamma_0)-\hat{u}_{t_0}(\gamma_0)=M_h-V_h(\gamma_0)>\max_{\Omega_h}V_h(z)-V_h(\gamma_0)\geq 0,$$
	the only case is  ${u}_h^{t_0}(\gamma_0)={u}_{t_0}(\gamma_0)=0.$ 
	Hence, there exists some $\lambda_0\neq t_0,$ such that
	\begin{multline*}
		M_h=\hat{u}_h^{t_0}(\gamma_0)-\hat{u}_{t_0}(\gamma_0)+V_h(\gamma_0)=\\=
		\sum_{l\neq t_0}\left({u}_{l}(\gamma_0)-{u}_h^{l}(\gamma_0) \right)+ V_h(\gamma_0)={u}_{\lambda_0}(\gamma_0)-\sum_{l\neq t_0}{u}_h^{l}(\gamma_0) + V_h(\gamma_0).
	\end{multline*}
	Recalling that $M_h=R_h$, we  write the following inequality
	$$
	{u}_{\lambda_0}(\gamma_0)-\sum_{l\neq t_0}{u}_h^{l}(\gamma_0) + V_h(\gamma_0)=M_h\geq \hat{u}_{\lambda_0}(\gamma_0)-\hat{u}_h^{\lambda_0}(\gamma_0)+V_h(\gamma_0),
	$$
	which in turn gives $2\sum\limits_{l\neq \lambda_0}{u}_h^l(\gamma_0)\leq 0,$ and therefore
	${u}_h^l(\gamma_0)=0,$ for all $l\neq \lambda_0.$ Hence,
	$$
	M_h={u}_{\lambda_0}(\gamma_0)-\sum_{l\neq t_0}{u}_h^{l}(\gamma_0) + V_h(\gamma_0)={u}_{\lambda_0}(\gamma_0)-{u}_h^{\lambda_0}(\gamma_0) + V_h(\gamma_0).
	$$
	This suggests us to apply the same approach as above and using the fact that $$L_h(\hat{u}_{\lambda_0}(\gamma_0)-\hat{u}_h^{\lambda_0}(\gamma_0) + V_h(\gamma_0))\geq 0,$$ we obtain
	\begin{equation}\label{nbr-cond2}
	M_h=\hat{u}_{\lambda_0}(\gamma_0)-\hat{u}_h^{\lambda_0}(\gamma_0)+ V_h(\gamma_0)=\hat{u}_{\lambda_0}(\theta)-\hat{u}_h^{\lambda_0}(\theta)+ V_h(\theta)>\max_{\Omega_h}V_h(z),
	\end{equation}
	 for all $\theta\in nbr(\gamma_0).$

	Thus, continuing this process all the time for the neighbor points, in view of \eqref{nbr-cond1} and \eqref{nbr-cond2}, we observe that for every mesh point $\gamma$ there always exists some $l_\gamma\in\overline{1,m}$ such that: 
	$$\mbox{either}\;\;
	\hat{u}_{l_\gamma}(\gamma)-\hat{u}_{h}^{l_\gamma}(\gamma)=M_h-V_h(\gamma)>0, \;\;\mbox{or}\;\; 
	\hat{u}_{l_\gamma}(\gamma)-\hat{u}_{h}^{l_\gamma}(\gamma)=V_h(\gamma)-M_h<0.
	$$
	On the other hand, it is clear that sooner or later, we will reach the boundary $\partial\Omega_h$ and this will  give a contradiction, because for every $\gamma\in\partial\Omega_h,$ and $l=\overline{1,m}$ we have $$\hat{u}_{l}(\gamma)-\hat{u}_{h}^{l}(\gamma)=\hat{u}_h^{l}(\gamma)-\hat{u}_{l}(\gamma)=0.$$ From this  we conclude that the only possibility is $$M_h\leq\max_{\Omega_h}V_h(z).$$
	In the light of Lemma \ref{lemma2}, the inequality $M_h\leq \max_{\Omega_h}V_h(z)$  implies  $R_h\leq \max_{\Omega_h}V_h(z)$ and vice-versa. Recalling the definition of $M_h$ and $R_h$  for arbitrary $l\in\overline{1,m},$ and $z\in\Omega_h$ we get
	\begin{align}
	\begin{cases}
    \hat{u}_l(z)-\hat{u}_h^l(z)+V_h(z)\leq \max_{\Omega_h}V_h(z),\\
    \hat{u}_h^l(z)-\hat{u}_l(z)+V_h(z)\leq \max_{\Omega_h}V_h(z).
	\end{cases}
	\end{align}
	This leads to 
	$$
	|\hat{u}_l(z)-\hat{u}_h^l(z)|\leq \max_{\Omega_h}V_h(z)-\min_{\Omega_h}V_h(z).
	$$
	In view of a function $V_h(z)$ we obtain 
	$$
	|\hat{u}_l(z)-\hat{u}_h^l(z)|\leq D_{\Omega}\cdot\sum_{l=1}^{m}||L_h u_l-\Delta u_l||_{L^\infty(\Omega)},
	$$
	for all $l=1,2,\dots,m,$ where $D_{\Omega}=\frac{a^2}{2}.$  
	This in turn implies that for every $z\in\Omega_h$ and $l=\overline{1,m}$ 
	we have
	$$
	|{u}_l(z)-{u}_h^l(z)|\leq 2D_{\Omega}\cdot\sum_{l=1}^{m}||L_h u_l-\Delta u_l||_{L^\infty(\Omega)}.
	$$
	
	Finally, we can write 
	$$
	||{u}_l-{u}_h^l||_{L^\infty(\Omega_h)}\leq 2D_{\Omega}\cdot\sum_{l=1}^{m}||L_h u_l-\Delta u_l||_{L^\infty(\Omega)},
	$$
	for every $l=1,2,\dots,m.$ This completes the proof. 
\end{proof}

\begin{corollary}
It is clear that due to Proposition \ref{scheme_conv}, we have  $u_h^l\to u_l,$ for every $l=1,2,\dots,m,$  whenever $h\to0,$ provided each component $u_l\in C^2(\Omega).$
\end{corollary}

\begin{corollary}
Assume $u_l\in C^4(\Omega),$ for all $l=1,2,\dots,m,$ then the Taylor expansion for the Laplacian operator yields  $L_hu_l-\Delta u_l=O(h^2).$ This together with Proposition \ref{scheme_conv} implies the following asymptotic decay: $$||{u}_l-{u}_h^l||_{L^\infty(\Omega_h)}=O(h^2).$$
\end{corollary}

 Similar convergence rates have been obtained in \cite{one-phase-error-izvestiya,one-phase-error,American-valuationl}  for the difference schemes of one-phase obstacle-like problems.

%\subsection*{Acknowledgment}

%\newpage
%\bibliographystyle{amsalpha}%
\bibliographystyle{acm}%
\bibliography{multi_obst}%\label{sec-Ref}

\begin{thebibliography}{10}

\bibitem{arakelyan2015finite}
{\sc Arakelyan, A.}
\newblock A finite difference method for two-phase parabolic obstacle-like
  problem.
\newblock {\em Armenian Journal of Mathematics 7}, 1 (2015), 32--49.

\bibitem{Av-Raf-2016}
{\sc Arakelyan, A., and Barkhudaryan, R.}
\newblock A numerical approach for a general class of the spatial segregation
  of reaction-diffusion systems arising in population dynamics.
\newblock {\em arXiv preprint arXiv:1603.03196\/} (2016).

\bibitem{one-phase-error-izvestiya}
{\sc Arakelyan, A., Barkhudaryan, R., and Poghosyan, M.}
\newblock An error estimate for the finite difference scheme for one-phase
  obstacle problem.
\newblock {\em Journal of Contemporary Mathematical Analysis 46}, 3 (2011),
  131--141.

\bibitem{ABP2014}
{\sc Arakelyan, A., Barkhudaryan, R., and Poghosyan, M.}
\newblock {Numerical Solution of The Two-Phase Obstacle Problem by Finite
  Difference Method}.
\newblock {\em Armenian Journal of Mathematics 7}, 2 (2015), 164--182.

\bibitem{Avet-Henrik}
{\sc Arakelyan, A., and Shahgholian, H.}
\newblock Multi-phase quadrature domains and a related minimization problem.
\newblock {\em Potential Analysis 45}, 1 (2016), 135--155.

\bibitem{MR2961456}
{\sc Arakelyan, A.~G., Barkhudaryan, R.~H., and Poghosyan, M.~P.}
\newblock Finite difference scheme for two-phase obstacle problem.
\newblock {\em Dokl. Nats. Akad. Nauk Armen. 111}, 3 (2011), 224--231.

\bibitem{MR2563520}
{\sc Bozorgnia, F.}
\newblock Numerical algorithm for spatial segregation of competitive systems.
\newblock {\em SIAM J. Sci. Comput. 31}, 5 (2009), 3946--3958.

\bibitem{Farid}
{\sc Bozorgnia, F.}
\newblock Numerical solutions of a two-phase membrane problem.
\newblock {\em Applied Numerical Mathematics 61}, 1 (2011), 92--107.

\bibitem{faridperturb}
{\sc Bozorgnia, F.}
\newblock Perturbation formula for the two-phase membrane problem.
\newblock {\em Advances in Difference Equations 2011}, 1 (2011), 1--16.

\bibitem{farid-optimal-partition}
{\sc Bozorgnia, F.}
\newblock Optimal partitions for first eigenvalues of the laplace operator.
\newblock {\em Numerical Methods for Partial Differential Equations 31}, 3
  (2015), 923--949.

\bibitem{Mywork}
{\sc Bozorgnia, F., and Arakelyan, A.}
\newblock Numerical algorithms for a variational problem of the spatial
  segregation of reaction--diffusion systems.
\newblock {\em Applied Mathematics and Computation 219}, 17 (2013), 8863--8875.

\bibitem{farid-multiQD}
{\sc Bozorgnia, F., and Bazarganzadeh, M.}
\newblock Numerical schemes for multi phase quadrature domains.
\newblock {\em International Journal of Numerical Analysis \& Modeling 11}, 4
  (2014), 726--744.

\bibitem{farid-fem}
{\sc Bozorgnia, F., and Valdman, J.}
\newblock A fem approximation of a two-phase obstacle problem and its a
  posteriori error estimate.
\newblock {\em arXiv preprint arXiv:1606.01020\/} (2016).

\bibitem{bucur-multi}
{\sc Bucur, D., and Velichkov, B.}
\newblock Multiphase shape optimization problems.
\newblock {\em SIAM Journal on Control and Optimization 52}, 6 (2014),
  3556--3591.

\bibitem{Aram-Caff}
{\sc Caffarelli, L., Karakhanyan, A., and Lin, F.-H.}
\newblock The geometry of solutions to a segregation problem for nondivergence
  systems.
\newblock {\em Journal of Fixed Point Theory and Applications 5}, 2 (2009),
  319--351.

\bibitem{nochetto-residual}
{\sc Chen, Z., and Nochetto, R.~H.}
\newblock Residual type a posteriori error estimates for elliptic obstacle
  problems.
\newblock {\em Numerische Mathematik 84}, 4 (2000), 527--548.

\bibitem{one-phase-error}
{\sc Cheng, X.-l., and Xue, L.}
\newblock On the error estimate of finite difference method for the obstacle
  problem.
\newblock {\em Applied Mathematics and Computation 183}, 1 (2006), 416--422.

\bibitem{MR2146353}
{\sc Conti, M., Terracini, S., and Verzini, G.}
\newblock Asymptotic estimates for the spatial segregation of competitive
  systems.
\newblock {\em Adv. Math. 195}, 2 (2005), 524--560.

\bibitem{MR2151234}
{\sc Conti, M., Terracini, S., and Verzini, G.}
\newblock A variational problem for the spatial segregation of
  reaction-diffusion systems.
\newblock {\em Indiana Univ. Math. J. 54}, 3 (2005), 779--815.

\bibitem{MR2283921}
{\sc Conti, M., Terracini, S., and Verzini, G.}
\newblock Uniqueness and least energy property for solutions to strongly
  competing systems.
\newblock {\em Interfaces Free Bound. 8}, 4 (2006), 437--446.

\bibitem{MR2300320}
{\sc Crooks, E. C.~M., Dancer, E.~N., and Hilhorst, D.}
\newblock Fast reaction limit and long time behavior for a
  competition-diffusion system with {D}irichlet boundary conditions.
\newblock {\em Discrete Contin. Dyn. Syst. Ser. B 8}, 1 (2007), 39--44
  (electronic).

\bibitem{MR1900331}
{\sc Dancer, E.~N., and Zhang, Z.}
\newblock Dynamics of {L}otka-{V}olterra competition systems with large
  interaction.
\newblock {\em J. Differential Equations 182}, 2 (2002), 470--489.

\bibitem{MR1459414}
{\sc Dancer, N.}
\newblock Competing species systems with diffusion and large interactions.
\newblock {\em Rend. Sem. Mat. Fis. Milano 65\/} (1995), 23--33.

\bibitem{multigrid}
{\sc Gr{\"a}ser, C., and Kornhuber, R.}
\newblock Multigrid methods for obstacle problems.
\newblock {\em Journal of Computational Mathematics\/} (2009), 1--44.

\bibitem{American-valuationl}
{\sc Hu, B., Liang, J., and Jiang, L.}
\newblock Optimal convergence rate of the explicit finite difference scheme for
  american option valuation.
\newblock {\em Journal of computational and applied mathematics 230}, 2 (2009),
  583--599.

\bibitem{oberman-siam}
{\sc Oberman, A.~M.}
\newblock Convergent difference schemes for degenerate elliptic and parabolic
  equations: Hamilton--jacobi equations and free boundary problems.
\newblock {\em SIAM Journal on Numerical Analysis 44}, 2 (2006), 879--895.

\bibitem{PSU2012}
{\sc Petrosyan, A., Shahgholian, H., and Ural'ceva, N.~N.}
\newblock {\em Regularity of free boundaries in obstacle-type problems},
  vol.~136.
\newblock American Mathematical Soc., 2012.

\bibitem{repin2015}
{\sc Repin, S., and Valdman, J.}
\newblock A posteriori error estimates for two-phase obstacle problem.
\newblock {\em Journal of Mathematical Sciences 207}, 2 (2015), 324--335.

\bibitem{rodrigues-book}
{\sc Rodrigues, J.-F.}
\newblock {\em Obstacle problems in mathematical physics}, vol.~134.
\newblock Elsevier, 1987.

\bibitem{MR2417905}
{\sc Squassina, M.}
\newblock On the long term spatial segregation for a competition-diffusion
  system.
\newblock {\em Asymptot. Anal. 57}, 1-2 (2008), 83--103.

\bibitem{Osher}
{\sc Tran, G., Schaeffer, H., Feldman, W.~M., and Osher, S.~J.}
\newblock An l\^{}1 penalty method for general obstacle problems.
\newblock {\em SIAM Journal on Applied Mathematics 75}, 4 (2015), 1424--1444.

\end{thebibliography}

\end{document}